\documentclass[12pt]{arXiv_class}
\usepackage{amsmath}
\usepackage[utf8]{inputenc}
\usepackage[T1]{fontenc}
\usepackage[utf8]{inputenc}
\usepackage{amsmath}
\usepackage{amsfonts}
\usepackage{amssymb}
\usepackage{xcolor}
\usepackage{stmaryrd}
\usepackage{tikz}
\usepackage{pgfplots}
\usepackage{hyperref} 

\usepackage{changes}
\definechangesauthor[name={joaquim garcia}, color=orange]{jg}
\definechangesauthor[name={mario souto}, color=red]{ms}
\setremarkmarkup{(#2)}

\usepackage{amssymb}

\renewcommand{\algorithmiccomment}[1]{\bgroup\hfill$\triangleright$~#1\egroup}

\newcommand*{\vertbar}{\rule[-1ex]{0.5pt}{2.5ex}}

\title{Exploiting Low-Rank Structure in Semidefinite Programming by Approximate Operator Splitting}
\author{Mario Souto, Joaquim D. Garcia and \'Alvaro Veiga}

\begin{document}
\maketitle

\begin{abstract}
In contrast to many other convex optimization classes, state-of-the-art semidefinite programming solvers are still unable to efficiently solve large scale instances. This work aims to reduce this scalability gap by proposing a novel proximal algorithm for solving general semidefinite programming problems. The proposed methodology, which is based on the primal-dual hybrid gradient method, allows for the presence of linear inequalities without the need to add extra slack variables and avoids solving a linear system at each iteration. More importantly, it simultaneously computes the dual variables associated with the linear constraints. The main contribution of this work is that it achieves a substantial speedup improvement by effectively adjusting the proposed algorithm in order to exploit the low-rank property inherent to several semidefinite programming problems. This modification is the key element that allows the operator splitting method to efficiently scale to larger instances. Convergence guarantees are presented along with an intuitive interpretation of the algorithm. Additionally, an open-source semidefinite programming solver called \texttt{ProxSDP} is made available and its implementation details are discussed. Case studies are presented in order to evaluate the performance of the proposed methodology.


\end{abstract}

\section{Introduction}

\subsection{Motivation and contributions}

Semidefinite programming (SDP) plays an important role in the field of convex optimization and subsumes several classes of optimization problems such as linear programming (LP), quadratic programming (QP) and second-order cone programming (SOCP). As a consequence, the range of applications to which SDP can be applied is wide and constantly expanding. In addition to being a general framework for convex problems, SDP is also a powerful tool for building tight convex relaxations of NP-hard problems. This property has significant practical consequences for the approximation of a wide range of combinatorial optimization problems and potentially to all constraint satisfaction problems \cite{raghavendra2008optimal}.


In practice, if one is interested in solving an SDP problem, it is crucial to have access to a fast, reliable and memory efficient software available. Unfortunately, in comparison to other convex optimization classes, the currently available SDP solvers are not as efficient as their counterparts. All these elements suggest that the development of an efficient algorithm and software for solving SDP problems would provide a noteworthy contribution. In this sense, the contributions of this paper are the following.
\begin{itemize}
\item A first order proximal algorithm for solving general SDP problems based on the \textit{primal-dual hybrid gradient} (PDHG) \cite{chambolle2011first} is proposed. The main advantage of this methodology, in comparison to other operator splitting techniques, is that it computes the optimal dual variables along with the optimal primal solution. Additionally, the algorithm does not require the solving of a linear system at every iteration and it allows for the presence of linear inequalities without the need to introduce additional variables into the problem.

\item Inspired by the approximate proximal point algorithm, a modified version of the PDHG that can exploit the low-rank property of SDP is proposed. For several problems of interest, this modification makes PDHG competitive with interior-point methods, in some cases providing a speed improvement of an order of magnitude. For problems with low-rank structure, the proposed algorithm is able to solve instances with dimensions that were still unattainable to interior-point methods in less than ten minutes, up to a $5,000 \times 5,000$ sized semidefinite matrix.

\item An open source SDP solver, called \texttt{ProxSDP}, is made publicly available. The goal of developing and providing this software is to both make the results of this paper reproducible and to foster the use of semidefinite programming in different fields.
\end{itemize} 	

The remainder of this section will cover some historical background on semidefinite programming, the current solution methodologies and introduce some notations. Section 2 will introduce the PDHG algorithm in the context of semidefinite programming. Section 3 will show how to modify the PDHG method in order to exploit the low-rank structure of the problem. In section 4, three case studies from different domains are considered in order to validate the proposed methodology.

\subsection{Notations}

In this work, we make use of the following notation. The symbol $X \succeq 0$ means that the matrix $X$ lies on the positive semidefinite cone, i.e. $y^* X y \geq 0  \hspace{0.2cm} \forall \hspace{0.1cm} y \in \mathbb{R}^n$. The symbol $\mathbb{S}^n$ represents the set of all $n \times n$ symmetric matrices and $\mathbb{S}_+^n$ is the set of all $n \times n$ symmetric positive semidefinite matrices. The symbol $\left\Vert \cdot \right\Vert_2$ denotes both the Euclidean norm for vectors and the \textit{spectral} norm for matrices. Additionally, $\left\Vert \cdot \right\Vert_F$ represents the Frobenius norm for matrices.

Given a function $f : \mathbb{R}^n \mapsto \mathbb{R} \cup \{\infty\}$, the associated subdifferential operator is defined as follows:
\begin{equation*}
\begin{aligned}
& \partial f = \{(x, g) : x \in \mathbb{R}^n, \hspace{0.1cm} f(y) \geq f(x) + g^T (y-x) \hspace{0.1cm} \forall \hspace{0.1cm} y \in \textbf{dom} \hspace{0.05cm} f \}.
\end{aligned}\label{eq_sdp}
\end{equation*}
The subdifferential operator evaluated at a point $x \in \mathbb{R}^n$ gives a set that is denoted by $\partial f (x)$, which is  called the subdifferential of $f$ at $x$ and is given by
\begin{equation*}
\begin{aligned}
& \partial f(x) = \{g : g^T (y-x) \leq f(y) - f(x) \hspace{0.1cm} \forall \hspace{0.1cm} y \in \textbf{dom} \hspace{0.05cm} f\}.
\end{aligned}\label{eq_sdp}
\end{equation*}
A subgradient of $f$ at $x$ is any point in the subdifferential of $f$ at $x$, i.e. $g \in \partial f(x)$. The inverse of the subdifferential operator, denoted by $(\partial f)^{-1}$, is defined as follows:
\begin{equation*}
\begin{aligned}
& (\partial f)^{-1} = \{(g, x) : (x, g) \in \partial f \}.
\end{aligned}\label{eq_sdp}
\end{equation*}

\subsection{Semidefinite programming formulation}

In this work we are going to consider a formulation of semidefinite programming where the inequalities are explicitly stated, thereby avoiding the use of slacks variables. The formulation referred to as \textit{general} SDP form is defined as the following
\begin{equation}
\begin{aligned}
& \underset{X \in \mathbb{S}^n}{\text{minimize}}
& & \textbf{tr}(CX)\\
& \text{subject to}
& & \mathcal{A}(X) = b,\\
&&& \mathcal{G}(X) \leq h,\\
&&& X \succeq 0.
\end{aligned}\label{eq_sdp_general}
\end{equation}
where the operators $\mathcal{A}: \mathbb{S}^n_+ \rightarrow \mathbb{R}^m$ and $\mathcal{G}: \mathbb{S}^n_+ \rightarrow \mathbb{R}^p$ are given by
\begin{equation*}
\begin{aligned}
&  \mathcal{A}(X) = \begin{bmatrix}
   \textbf{tr}(A_1 X)  \\
    \textbf{tr}(A_2 X) \\
    \vdots \\
    \textbf{tr}(A_m X) 
\end{bmatrix}, \hspace{0.3cm}
\mathcal{G}(X) = \begin{bmatrix}
   \textbf{tr}(G_1 X)  \\
    \textbf{tr}(G_2 X) \\
    \vdots \\
    \textbf{tr}(G_p X)
\end{bmatrix}
\end{aligned}
\end{equation*}
and the problem data are the symmetric matrices $A_1, \dots, A_m, G_1, \dots, G_p, C \in \mathbb{S}^n$ and the vectors $b \in \mathbb{R}^m$ and $h \in \mathbb{R}^p$. In this semidefinite programming formulation, one wants to minimize a linear function subjected to a set of $m$ linear equality constraints and $p$ linear inequalities constraints, where the decision variable is an $n \times n$ symmetric matrix constrained to be on the positive semidefinite (p.s.d.) cone. 

The dual problem takes the form
\begin{equation*}
\begin{aligned}
& \underset{y \in \mathbb{R}^{m + p}}{\text{maximize}}
& & [b^T \hspace{0.1cm} h^T] \hspace{0.1cm} y\\
& \text{subject to}
& & \mathcal{A}^*(y) + \mathcal{G}^*(y) \preceq C,\\
&&& y_j \leq 0, \hspace{0.2cm} \forall \hspace{0.1cm} j = m+1, \dots, p.
\end{aligned}\label{eq_sdp}
\end{equation*}
where the conjugate operators $\mathcal{A}^*: \mathbb{R}^m \rightarrow \mathbb{S}^n_+$ and $\mathcal{G}^*: \mathbb{R}^p \rightarrow \mathbb{S}^n_+$ are given by $\mathcal{A}^*(y) = \sum_{i=1}^{m} y_{i} A_i$ and $\mathcal{G}^*(y) =\sum_{j=1}^{p} y_{j+m} G_j$ respectively. Despite being very similar to linear programming, strong duality does not always hold for SDP. For a comprehensive analysis on the duality of semidefinite programming the reader should refer to the work of Boyd and Vandenberghe \cite{vandenberghe1996semidefinite}.

\subsection{Applications}

In addition to the theoretical motivation, several problems of practical interest lie precisely in the SDP class. As more applications are found, an efficient method for solving large scale SDP problems is required. This section, briefly covers some representative applications that have been successfully solved by SDP. For a more complete list of applications and SDP problems, the reader can refer to \cite{wolkowicz2012handbook, vandenberghe1999applications}.

Traditionally, semidefinite programming has been widely used in control theory. Classic applications such as the stability of dynamic systems and stochastic control problems \cite{lur1957some} have motivated the development of SDP over decades. Modern applications such as motion for humanoid robots \cite{kuindersma2016optimization} use SDP as a core element for control. In several cases, SDP problems in control can be formulated in the form of \textit{linear matrix inequalities} (LMI) \cite{bellman1963systems}. In this sense, for some particular applications in control, there is no need to use general SDP algorithms since some LMI problems have closed form solutions, although adding a little complexity to the problem might render the closed form solutions useless. As a consequence, interior-point methods that were particularly developed to solve LMI problems were proposed \cite{vandenberghe2005interior}. Extensive literatures on LMI can be found in \cite{boyd1994linear, scherer2000linear}. 

Subsequently, several fields of study have found SDP to be a powerful tool for solving complex problems. For instance, the use of SDP in power systems allowed for deriving tight bounds and solutions for more realistic optimal power flow models with alternating current networks \cite{lavaei2012zero}. In chip design,  transistor sizing was also optimized with the use of SDP \cite{vandenberghe1998optimizing, vandenberghe1997optimal}. In the field of structural truss layout, the use of SDP has been popularized after the seminal work of Ben-Tal and Nemirovski \cite{ben1997robust}. This latter work also presented the first ideas that subsequently lead to the expanding field of optimization under uncertainty \cite{ben1998robust}, where SDP is also used to approximate chance constraints \cite{zymler2013distributionally}.

A remarkable property of SDP is the ability to build tight convex relaxations to NP-hard problems. This technique, which is also known as semidefinite relaxation (SDR), has been a powerful tool bridging convex and combinatorial optimization. In the early nineties, Lov\'asz and Schrijver developed a SDR for optimization problems with the presence of boolean variables \cite{lovasz1991cones}. Soon after, SDR gained momentum after the celebrated Goemans-Williamson randomized rounding method for the max-cut problem \cite{goemans1995improved}. Eventually, similar SDR approaches were proposed for other combinatorial problems, such as the max-3-sat \cite{karloff19977} and the traveling salesman problem \cite{cvetkovic1999semidefinite}. Recently, Cand\`es et al. proposed an SDR approach to the phase retrieval problem \cite{candes2015phase}. Such relaxations generally square the original number of decision variables through a technique called \textit{lifting} \cite{lovasz2003semidefinite}. From the algorithmic perspective, this can quickly make moderate size instances computationally challenging.

In the last decade, applications in machine learning have challenged traditional SDP methods at solving remarkably large-scale instances \cite{de2007deploying}. The problem size $n$ is usually associated with the sample size, which can easily reach millions. One of the most celebrated applications is the \textit{matrix completion} problem \cite{candes2012exact} which became very popular with the \textit{Netflix prize} \cite{bennett2007netflix}. In graphical models, the problem of covariance selection, which is a powerful tool for modeling dependences between random variables, can also be formulated as an SDP problem \cite{wainwright2006log}. In statistical learning, finding the best model that combines different positive definite kernels, which is also known as kernel learning, can be achieved by solving an SDP problem \cite{lanckriet2004learning, lanckriet2004statistical}.

For \textit{constraint satisfaction problems} (CSP), problems where one tries to satisfy as many constraints as possible, semidefinite programming also plays an essential role. In the work of Raghavendra \cite{raghavendra2008optimal}, it was shown that if the \textit{Unique Games} conjecture \cite{khot2002power} is true, then semidefinite programming achieves the best approximation for every CSP. Even though there is no consensus on the legitimacy of the Unique Games conjecture, the recent work by Khot et al. \cite{khot2018pseudorandom, dinur2018towards} provides new results suggesting the veracity of the conjecture. Those recent developments bring semidefinite programming back into the spotlight and imply that this class of convex optimization algorithms does have singular properties worth exploiting.

\subsection{SDP solution methods}

In the early days, general SDP problems were solved by the ellipsoid method \cite{shor1977cut, iudin1977informational} and subsequently by bundle methods \cite{hiriart2013convex}. After the advent of the first polynomial time interior-point algorithm for linear programming by Karmarkar et al. \cite{karmarkar1984new, adler1989implementation}, Nesterov and Nemirovski extended the interior-point methods for other classes of convex optimization problems \cite{nesterov1994interior}. Shortly afterwards, a range of interior-point methods to solve SDPs were proposed \cite{alizadeh1992optimization, helmberg1996interior, alizadeh1998primal}. The solvers CSDP \cite{borchers1999csdp} and MOSEK \cite{mosek2010mosek} use state-of-the-art interior points methods for solving general SDP problems. Up to medium size problems, this class of methods is preferable due to the fast convergence and high precision. However, as is inherent to all second-order methods, the use of interior-point algorithms may be prohibitive for solving large-scale instances. The main bottleneck is due to the cumbersome effort for computing and storing the Hessian at each iteration.

More recently, first-order methods have been widely used for applications in machine learning and signal processing. Even though first-order methods generally have slower convergence rates, their costs per iteration are usually small and they require less memory allocation \cite{boyd2011distributed}. These characteristics make first-order methods very appealing for large scale problems and, consequently, it has been an intense area of research in several fields, such as image processing \cite{heide2016proximal}. A great example of the use of first order methods is the conic solver SCS developed by O'Donoghue et al. \cite{o2016conic}, which can efficiently solve general conic optimization problems to modest accuracy. More recently, the use of the \textit{alternating direction method of multipliers} (ADMM) for specifically solving SDP has been proposed by \cite{madani2015admm}. 

For most of the algorithms, exploiting sparsity patterns in the decision variables is not as straightforward as it is for other classes of convex optimization problems. In this sense, a major recent contribution has been made by showing that sparsity can be exploited by means of chordal decomposition techniques \cite{fukuda2001exploiting, nakata2003exploiting, vandenberghe2015chordal, pakazad2018distributed}. This approach has enabled parallel implementations that can solve larger instances with the use of supercomputers \cite{fujisawa2012high}. In a series of works, Zhang and Lavaei presented SDP algorithms that can properly take advantage of the problem's sparsity  \cite{zhang2017sparse, zhang2017modified}. 









\section{A primal-dual operator splitting for SDP}\label{sec_fixed_point}

In this section, the proposed operator splitting method will be built from the first-order optimality conditions of the SDP in its general form (\ref{eq_sdp_general}). The strategy adopted to derive the algorithm translates the problem of finding a solution that satisfies the optimality conditions into a problem of finding a fixed point of a related monotone operator. This approach has been previously adopted with the purpose of designing new algorithms and developing alternative proofs for existing ones \cite{ryu2016primer}. Further information on the use of monotone operators in the context of convex optimization can be found in \cite{eckstein1989splitting, combettes2004solving, combettes2005signal, combettes2011proximal}.

Consider the general SDP form (\ref{eq_sdp_general}) where the problem constraints are encoded by indicator functions
\begin{equation}
\begin{aligned}
& \underset{X \in \mathbb{S}^n}{\text{minimize}}
& & \textbf{tr}(C X) + I_{\mathbb{S}^n_+}(X) + I_{{=b \atop \leq h}}(\mathcal{M}(X)), \hspace{0.3cm} \mathcal{M} = \begin{bmatrix} \mathcal{A} \\ \mathcal{G}\end{bmatrix}.
\end{aligned}\label{eq_sdp}
\end{equation}
Where the indicator functions are defined as the following:
\begin{equation*}
\begin{aligned}
&I_{\mathbb{S}^n_+}(X) = \left\{
                \begin{array}{ll}
                  0, \hspace{0.3cm} \text{if} \hspace{0.1cm} X \succeq 0, \\
                  \infty, \hspace{0.15cm} \text{otherwise,}
                \end{array}
              \right.
\end{aligned}
\end{equation*}
encodes the positive semidefinite cone constraint and 
\begin{equation*}
\begin{aligned}
& I_{{=b \atop \leq h}}(u) = I_{=b}(u_1) + I_{\leq h}(u_2), \hspace{0.5cm} u = [u_1 \hspace{0.2cm} u_2]^T,
\end{aligned}
\end{equation*}
encodes the linear constraints right-hand side for any $u \in \mathbb{R}^{m + p}$ such that
\begin{equation*}
\begin{aligned}
& I_{=b}(u_1) = \left\{
                \begin{array}{ll}
                  0, \hspace{0.3cm} \text{if} \hspace{0.1cm} u_1 = b, \\
                  \infty, \hspace{0.15cm} \text{otherwise,}
                \end{array}
              \right. \hspace{1cm} 
I_{\leq h}(u_2) = \left\{
                \begin{array}{ll}
                  0, \hspace{0.3cm} \text{if} \hspace{0.1cm} u_2 \leq h, \\
                  \infty, \hspace{0.15cm} \text{otherwise.}
                \end{array}
              \right.
\end{aligned}
\end{equation*}
for any $u_1 \in \mathbb{R}^{m}$ and $u_2 \in \mathbb{R}^{p}$.

\subsection{Optimality condition}

The first order optimality condition (\ref{eq_opt_system}) for the optimization problem (\ref{eq_sdp}) can be expressed as follows
\begin{equation}
	\begin{aligned}
	& 0 \in \partial \hspace{0.1cm} \textbf{tr}(C X) + \partial \hspace{0.1cm} I_{\mathbb{S}^n_+}(X) + \mathcal{M}^* (\partial \hspace{0.1cm} I_{{=b \atop \leq h}}(\mathcal{M}(X))).
	\end{aligned}\label{eq_opt_system}
\end{equation} 
By introducing an auxiliary variable $y \in \mathbb{R}^{m + p}$, the optimality condition can be recast as the following system of inclusions
\begin{equation*}
	\begin{aligned}
	& 0 \in \partial \hspace{0.1cm} \textbf{tr}(C X) + \partial \hspace{0.1cm} I_{\mathbb{S}^n_+}(X) + \mathcal{M}^*(y),\\
    & y \in \partial \hspace{0.1cm} I_{{=b \atop \leq h}} (\mathcal{M}(X)).
	\end{aligned}
\end{equation*}
By definition, the auxiliary variable $y$ represents the dual variable associated with the problem constraints. This statement is easily verifiable since $y \in \partial \hspace{0.1cm} I_{{=b \atop \leq h}}(\mathcal{M}(X))$, i.e. $y$ is a subgradient of $I_{{=b \atop \leq h}}(\mathcal{M}(X))$ at $X$. Since problem (\ref{eq_sdp}) is convex, finding a pair $(X^*, y^*)$ satisfying (\ref{eq_opt_system}) is equivalent to finding an optimal primal-dual pair for (\ref{eq_sdp}) as long as strong duality holds \cite{bauschke2011convex}.

Using the fact that $(\partial f)^{-1} = \partial f^*$, for an $f$ that is a convex closed proper \cite{rockafellar2015convex}, one can manipulate the second inclusion as follows
\begin{equation*}
	\begin{aligned}
    & y \in \partial \hspace{0.1cm} I_{{=b \atop \leq h}} (\mathcal{M}(X)) \iff \partial \hspace{0.1cm} I^*_{{=b \atop \leq h}} (y) \ni \mathcal{M}(X),\\
    & \hspace{2.85cm} \iff 0 \in \partial \hspace{0.1cm} I^*_{{=b \atop \leq h}} (y) - \mathcal{M}(X).
	\end{aligned}
\end{equation*}
Using this new expression, the system (\ref{eq_opt_system}) can be recast as 
$0 \in F(X, y)$, where $F$ is given by the following monotone operator
\begin{equation}
	\begin{aligned} F(X, y) = 
	\big( \partial \hspace{0.1cm} \textbf{tr}(C X) + \partial \hspace{0.1cm} I_{\mathbb{S}^n_+}(X) + \mathcal{M}^*(y) \hspace{0.05cm}, \hspace{0.1cm}  \partial \hspace{0.1cm} I_{{=b \atop \leq h}}^*(y) - \mathcal{M}(X) \big).
	\end{aligned}\label{eq_F}
\end{equation}
One can verify that $F$ is a monotone operator by noticing that $F$ is the sum of a subdifferential operator and a monotone affine operator \cite{bauschke2017convex}


This formulation of the inclusion (\ref{eq_opt_system}) implies that finding a zero of the underlying monotone operator $F$ is equivalent to finding an optimal primal-dual pair for the semidefinite programming problem (\ref{eq_sdp}). In the remainder of this section, a method for finding a zero for the operator $F$ will be established.

\subsection{Fixed point iteration}

Finding a zero of the monotone operator $F$ can be translated into finding a fixed point for the system $P(X, u) \in \alpha F(X, u) + P(X, u)$, where $P$ is a positive definite operator \cite{bauschke2011convex}. This formulation induces the following fixed point iteration
\begin{equation}
	\begin{aligned}
    & \big(X^{k}, u^{k}\big) \leftarrow \big(P + \alpha F\big)^{-1} P(X^{k-1}, u^{k-1}).
	\end{aligned}\label{eq_fixed_point}
\end{equation}
This iterative process is called the \textit{generalized proximal point method} and it is guaranteed to converge to a zero of $F$ if the matrix $P$ is positive definite and a fixed point for (\ref{eq_fixed_point}) exists \cite{ryu2016primer}. 

By choosing $P$ as
\begin{equation}
	\begin{aligned}
    & P = \begin{bmatrix}
           I & - \alpha \mathcal{M}^* \\[0.3cm]  - \alpha \mathcal{M} & I
        \end{bmatrix},
	\end{aligned}\label{eq_P}
\end{equation}
the fixed point inclusions can be expressed as 
\begin{equation*}
	\begin{aligned}
		& \left(X^{k-1} - \alpha \mathcal{M}^*(y^{k-1}), \hspace{0.1cm} y^{k-1} - \alpha \mathcal{M}(X^{k-1})\right) \in \alpha F(X^k, y^k) \hspace{-0.05cm} + \left(X^k - \alpha \mathcal{M}^*(y^k), \hspace{0.1cm} y^k - \alpha \mathcal{M}(X^k)\right)
	\end{aligned}
\end{equation*}
where further manipulation leads to the system
\begin{equation*}
	\begin{aligned}
		& \hspace{-0.1cm} \big( X^{k-1} - \alpha \mathcal{M}^*(y^{k-1}), \hspace{0.1cm} y^{k-1} + \alpha \mathcal{M}(2X^{k} - X^{k-1}) \big) \in \alpha \big(\partial \hspace{0.1cm} \textbf{tr}(C X^{k}) + \partial \hspace{0.1cm} I_{\mathbb{S}^n_+}(X^{k}), \hspace{0.1cm}\partial \hspace{0.1cm} I^*_{{=b \atop \leq h}} (y^{k}) \big) \hspace{-0.05cm} + \big(X^k, \hspace{0.1cm} y^k \big)
	\end{aligned}
\end{equation*}
which induces the following fixed point iteration:
\begin{equation*}
	\begin{aligned}
		& X^{k} \leftarrow \big(I + \alpha \partial \hspace{0.1cm} (\textbf{tr}(C \cdot) + I_{\mathbb{S}^n_+}))^{-1} (X^{k-1} - \alpha \mathcal{M}^{*}(y^{k -1}) \big) \\[0.1cm]
        & y^{k} \leftarrow \big(I + \alpha \partial \hspace{0.1cm} I^*_{{=b \atop \leq h}})^{-1} (y^{k-1} + \alpha \mathcal{M}(2X^{k} - X^{k-1}) \big).
	\end{aligned}\label{eq_cp_fixed_point}
\end{equation*}

This scheme is a particular case of the primal-dual hybrid gradient (PDHG)  proposed by Chambolle and Pock \cite{chambolle2011first, pock2009algorithm} which has been successfully applied to a wide range of image processing problems, such as image denoising and deconvolution \cite{sidky2012convex, vaiter2013robust, heide2016proximal}. Convergence is guaranteed as long as a solution to (\ref{eq_sdp}) exists and $0 < \alpha < 1 / \left\Vert \mathcal{M} \right\Vert_2$. The progress of the algorithm can be measured by the primal, dual and combined residuals as
\begin{equation*}
	\begin{aligned}
	& \epsilon_{\text{primal}}^k = \left\Vert \tfrac{1}{\alpha}(X^{k} - X^{k - 1}) - \mathcal{M}^*(y^k - y^{k - 1}) \right\Vert_F,\\
    & \epsilon_{\text{dual}}^k = \left\Vert \tfrac{1}{\alpha} (y^{k} - y^{k - 1}) - \mathcal{M}(X^k - X^{k - 1}) \right\Vert_2,\\
    & \epsilon_{\text{comb}}^k = \epsilon_{\text{primal}}^k + \epsilon_{\text{dual}}^k.
	\end{aligned}\label{eq_comb_residual}
\end{equation*}

\subsection{Resolvents and proximal operators in SDP}\label{subsec_prox}

To employ the fixed point iteration (\ref{eq_fixed_point}) one needs to compute both resolvent operators $(I + \alpha \partial \hspace{0.1cm} (\textbf{tr}(C \cdot) + I_{\mathbb{S}^n_+}))^{-1}$ and $(I + \alpha \partial \hspace{0.1cm} I^*_{{=b \atop \leq h}})^{-1}$. For any  convex function $f:\mathbb{R}^{m \times n} \rightarrow \mathbb{R} \cup \{+ \infty\}$ and $\alpha > 0$, the resolvent operator associated with the subdifferential operator $\partial  f$ is given by
\begin{equation*}
	\begin{aligned}
	& z = (I + \alpha \partial f)^{-1} (v) \iff 0 \in  \partial f (z) + \tfrac{1}{\alpha} (z-v) \iff z = \underset{x}{\text{argmin}}\Big\{f(x) + \tfrac{1}{2 \alpha}||x - v||_{\text{2}}^2\Big\}.
	\end{aligned}
\end{equation*}
Therefore, for a general convex function $f$, the resolvent of $\partial f$ can be expressed as the solution of an associated convex optimization problem. This mapping is usually referred as the  \textit{proximal operator} of $f$. The constant
 $\alpha > 0$ is a parameter that controls the trade-off between moving towards the minimizer of $f$ and shifting in the direction of $v$.


Methods based on proximal operators, such as the proximal gradient descent \cite{combettes2011proximal}, are particularly interesting when the proximal operator has a known closed-form. When the objective function is not differentiable, proximal algorithms are generally a good alternative to subgradient based methods, such as ISTA for linear inverse problems \cite{beck2009fast}. For a deeper review of proximal algorithms and more details on the resolvent calculus employed in the forthcoming sections the reader should refer to \cite{parikh2014proximal}. In the following, the proximal operators associated with (\ref{eq_fixed_point}) are going to be analyzed in more detail.

\subsubsection{Box constraints}

The resolvent associated with $\partial I_{{=b \atop \leq h}}$ is simply given by the projection onto the box constraints 
\begin{equation*}
	\begin{aligned}
	& \big(I + \alpha \partial I_{{=b \atop \leq h}}\big)^{-1} (u) = \textbf{proj}_{{=b \atop \leq h}}(u) = \begin{bmatrix}
           \textbf{proj}_{=b}(u_1) \\[0.3cm]  \textbf{proj}_{\leq h}(u_2)
        \end{bmatrix}
= \begin{bmatrix}
           b \\[0.3cm]  \min{\{u_2,h\}}
        \end{bmatrix},
	\end{aligned}
\end{equation*}
where $u_1 \in \mathbb{R}^p$ and $u_2 \in \mathbb{R}^m$ and $\min{\{u_2,h\}}$ is the point-wise minimum. Additionally, for a convex function $f$ , the \textit{extended Moreau decomposition} \cite{moreau1962decomposition} gives the identity
\begin{equation*}
	\begin{aligned}
	& u = \big(I + \alpha \partial f^*\big)^{-1} (u) + \alpha \hspace{0.05cm} \big(I + \partial f / \alpha \big)^{-1} (u / \alpha).
	\end{aligned}
\end{equation*}
Therefore one concludes that
\begin{equation}
	\begin{aligned}
	& \big(I + \alpha \partial I_{{=b \atop \leq h}}^*\big)^{-1} (u) = u - \alpha \hspace{0.05cm} \textbf{proj}_{{=b \atop \leq h}}(u / \alpha).
	\end{aligned}\label{eq_resolvent_1}
\end{equation}

\subsubsection{Positive semidefinite cone}

Similarly, the resolvent associated with the positive semidefinite constraint is given by the Euclidean projection onto the positive semidefinite cone. Let $S \in \mathbb{S}^n$, the projection onto the set $\{X : X \succeq 0 \}$ has the closed form
\begin{equation*}
	\begin{aligned}
	& \big(I + \alpha \partial I_{\mathbb{S}_+^n}\big)^{-1} (S) = \textbf{proj}_{\mathbb{S}_+^n}(S) = \sum_{i=1}^n \text{max}\{0, \lambda_i\} u_i u_i^{T},
	\end{aligned}
\end{equation*}
where $S = \sum_{i=1}^n \lambda_i u_i u_i^{T}$ is the eigenvalue decomposition of the symmetric matrix $S$.

\subsubsection{Trace}

Given the symmetric matrices $C$ and $S$, the resolvent associated with the trace function is given by the formula
\begin{equation*}
	\begin{aligned}
	& \big(I + \alpha \partial \textbf{tr}(C \cdot)\big)^{-1} (S) = S - \alpha C.
	\end{aligned}
\end{equation*}
Unlike the majority of cases, the resolvent associated with the trace function plus any convex function $g$ is given by the left composition as in
\begin{equation*}
	\begin{aligned}
	& \big(I + \alpha \partial (g + \textbf{tr}(C \cdot)\big)^{-1} (S) = \big(I + \alpha \partial g \big)^{-1} \circ \big( I + \alpha \partial \textbf{tr}(C \cdot)\big)^{-1} (S) = \big(I + \alpha \partial g \big)^{-1} \big(S - \alpha C \big).
	\end{aligned}
\end{equation*}
Consequently,
\begin{equation}
	\begin{aligned}
& \big(I + \alpha \partial \hspace{0.1cm} (\textbf{tr}(C \cdot) + I_{\mathbb{S}^n_+})\big)^{-1} (S) = \textbf{proj}_{\mathbb{S}_+^n} (S - \alpha C).
	\end{aligned}\label{eq_resolvent_2}
\end{equation}

\subsection{PD-SDP}

Algorithm 1, which is referred to as \texttt{PD-SDP} or \texttt{P}rimal-\texttt{D}ual \texttt{S}emi\texttt{D}efinite \texttt{P}rogramming, matches the fixed point iteration (\ref{eq_fixed_point}) and the resolvents in its closed forms (\ref{eq_resolvent_1}, \ref{eq_resolvent_2}).  In this particular setting, the primal-dual method turns out to be a very simple routine. As illustrated in Algorithm 1, the method avoids explicitly solving a linear system or a convex optimization problem at each iteration. One only needs a subroutine to evaluate the resolvents (\ref{eq_resolvent_1}, \ref{eq_resolvent_2}) and access the abstract linear operator for $\mathcal{M}$ and its adjoint.

\def\smath#1{\text{\scalebox{.9}{$#1$}}}

\begin{algorithm}[H]
	\caption{\texttt{\texttt{PD-SDP}}}
	\label{Algorithm} 
	\begin{algorithmic}
		\STATE{\textbf{Given:} $\mathcal{M}$, $b \in \mathbb{R}^m$, $h \in \mathbb{R}^p$ and $C \in \mathbb{S}^n$. \vspace{0.2cm}}
	    \WHILE{$\epsilon_{\text{comb}}^k > \epsilon_{\text{tol}}$ \vspace{0.2cm}}
        \STATE{$X^{k + \smath{1}} \hspace{0.1cm} \leftarrow 		\textbf{proj}_{\mathbb{S}_+^n} (X^k - \alpha ( \mathcal{M}^* (y^k) + C))$ \vspace{0.2cm}}\COMMENT{Primal step}
        \STATE{$y^{k + \smath{1/2}} \hspace{-0.05cm} \leftarrow y^{k} + \alpha \mathcal{M}(2 X^{k + \smath{1}} - X^k)$ \vspace{0.1cm}}\COMMENT{Dual step part 1}
        \STATE{$y^{k + \smath{1}} \hspace{0.2cm} \leftarrow y^{k + \smath{1/2}} - \alpha \hspace{0.05cm} \textbf{proj}_{{=b \atop \leq h}}(y^{k + \smath{1/2}} / \alpha)$}\COMMENT{Dual step part 2 \vspace{0.2cm}}
        \ENDWHILE{\vspace{0.2cm}}
		\RETURN{$\big(X^{k + 1}, y^{k + 1}\big)$}
	\end{algorithmic}
\end{algorithm}

In \cite{chambolle2011first}, it was shown that PDHG achieves a $\mathcal{O}(1/k)$ convergence rate for non-smooth problems, where $k$ is the number of iterations. Similar convergence rates can be achieved by other operator splitting methods such as Tseng's ADM \cite{tseng1991applications} or ADMM \cite{glowinski1975approximation, eckstein1992douglas}. However, the \texttt{PD-SDP} has the advantage of offering the optimal dual variable as a by-product of the algorithm.

The computational complexity of each loop is dominated by the projection onto the positive semidefinite cone. In the most naive implementation, each iteration will cost $\mathcal{O}(n^3)$ operations, where $n$ is the dimension of the p.s.d. matrix. If one knew the number of positive eigenvalues $r$, at each iteration, the computational cost could be reduced to $\mathcal{O}(n^2 r)$. In this work we are going to refer to $r$ as the \textit{target-rank} of a particular iteration. Unfortunately, in practice, one does not have access to the \textit{target-rank}. However, as will be shown in the next section, there is no need to know the \textit{target-rank} in advance. Even more surprisingly, faster running times can be achieved by underestimating the target rank to some extent. 

\section{Speeding up with inexact solves}

So far, we have proposed a first-order method for solving general SDP problems. However, the projection onto the positive semidefinite cone is an obstacle to make the algorithm scalable for larger instances. In this section, we are going to explore how to take advantage of a low-rank structure, even if the target rank is unknown. 

\subsection{Low-rank approximation}

It is well-known that SDP solutions very often exhibit a low-rank structure. More precisely, as shown by Barvinok \cite{barvinok1995problems} and Pataki \cite{pataki1998rank}, any SDP with $m$ equality constraints has an optimal solution with a rank of at most $\sqrt{2m}$. In practice, for several SDP problems, it is frequently observed that the optimal solution has an even smaller rank. This phenomenon is notably present in SDPs generated by a semidefinite relaxation, where the solution ideally has low rank. In several cases, even if the relaxation is inexact, the rank of the solution is usually substantially small.

This property has motivated a series of nonconvex methods aiming to exploit the low-rank structure of the problem \cite{zhang2011penalty, zhao2012approximation, yuan2016proximal}. For instance, one can encode the positive semidefinite constraint as a matrix factorization of the type $X=V^{T} V$ where $V \in \mathbb{R}^{r \times n}$ and $r$ is the target rank. This technique was proposed a decade ago by Burer and Monteiro \cite{burer2003nonlinear} and since then it has been one of the main tools for tackling the scalability of low-rank SDPs. This matrix factorization approach has been successfully applied to large-scale computer vision \cite{shah2016biconvex} and combinatorial optimization problems \cite{wang2017mixing}. Unfortunately, by resorting to this approach, one loses convexity and all the associated guarantees.

The main bottleneck of \texttt{PD-SDP} and any other convex optimization methods for solving SDPs is computing the eigenvalue decomposition. A natural approach to overcome this issue is to make use of low-rank matrix approximation techniques in place of computing the full matrix decomposition. Recent work by Udell, Tropp et al. \cite{yurtsever2017sketchy} uses matrix sketching methods \cite{halko2011finding, tropp2017practical} to successfully find approximate solutions to low-rank convex problems. While their methodology possesses several  advantages, such as optimal storage, it does require all solutions to be low-rank in order to guarantee convergence to an optimal solution. In contrast, the methodology proposed in this paper exploits the low-rank structure of the problem whenever possible, but it also converges to an optimal solution even in the presence of full-rank solutions.

As was showed by Eckart and Young \cite{eckart1936approximation}, the best rank-$r$ approximation of symmetric matrices, for both the Frobenius and the spectral norms, is given by the truncated eigenvalue decomposition. Inspired by this result, the \textit{approximate} projection onto the positive semidefinite cone is given by
\begin{equation}
	\begin{aligned}
	& \textbf{aproj}_{\mathbb{S}^n_+}(X, r) = \sum_{i=1}^r \text{max}\{0, \lambda_i\} u_i u_i^{T},
	\end{aligned}\label{eq_trunc_eig}
\end{equation}
where $X$ is a symmetric matrix, $r$ is its \textit{target-rank} and $\lambda_1 \geq \cdots \geq \lambda_r$ are the eigenvalues with the $r$ largest real values. It is important to notice that, despite being different from the Euclidean projection, $\textbf{aproj}_{\mathbb{S}^n_+}(X, r)$ does project the matrix $X$ onto the p.s.d. cone. In other words, the truncated projection maps onto the p.s.d. cone but not necessarily onto the closest point, according to the Frobenius norm, as illustrated in Figure 1. 
\begin{figure}[H]
\centering
\includegraphics[width=0.33\linewidth]{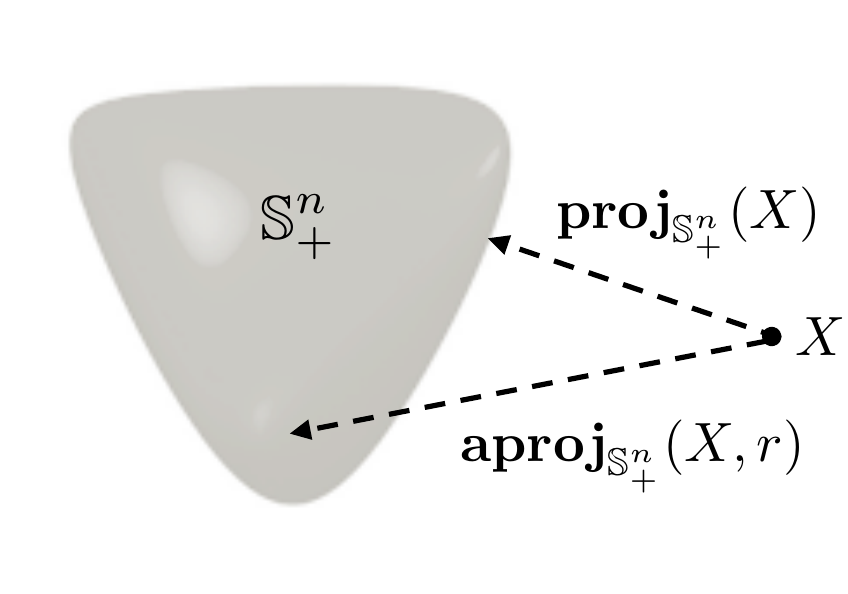}
\caption{Comparison of Euclidean projection onto the positive semidefinite cone, denoted by $\textbf{proj}_{\mathbb{S}^n_+}(X)$, and the truncated projection given by $\textbf{aproj}_{\mathbb{S}^n_+}(X, r)$.}
\end{figure}

If the target-rank $r$ equals the number of nonzero eigenvalues, both the truncated and the full projection will be equivalent. Otherwise, if the target-rank $r$ is smaller than the number of nonzero eigenvalues, the truncated projection will be only an approximation of the exact projection.  In this case, according to the Eckart–Young–Mirsky theorem \cite{eckart1936approximation}, the approximation error can be expressed as the sum of the eigenvalues that were left out by the truncated projection as
\begin{equation*}
	\begin{aligned}
	& \left\Vert \textbf{proj}_{\mathbb{S}^n_+}(X) - \textbf{aproj}_{\mathbb{S}^n_+}(X, r) \right\Vert_F^2 = \sum_{i=r + 1}^{n} \max\{\lambda_i, 0\}.
	\end{aligned}
\end{equation*}
For practical purposes, the approximation error can be bounded in terms of the smaller eigenvalue computed by the truncated projection as the following
\begin{equation}
	\begin{aligned}
	& \left\Vert \textbf{proj}_{\mathbb{S}^n_+}(X) - \textbf{aproj}_{\mathbb{S}^n_+}(X, r) \right\Vert_F^2 \leq (n - r) \max\{\lambda_r, 0\}.
	\end{aligned}\label{eq_error_approx}
\end{equation}

The partial eigenvalue decomposition  (\ref{eq_trunc_eig}) can be efficiently computed via power iteration algorithms or Krylov subspace methods \cite{golub2012matrix, higham1988matrix}. Computational routines are freely available in almost every programming language for numerical computing \cite{lehoucq1998arpack, anderson1999lapack}.


\subsection{Convergence checking and target rank update}

As previously noted, the \texttt{PD-SDP} method can be seen as a fixed point iteration of a monotone operator \cite{he2010convergence}. In this sense, replacing the exact projection onto the positive semidefinite cone by its approximation (\ref{eq_trunc_eig}) will result in an inexact iteration as the following
\begin{equation}
	\begin{aligned}
    & \big(X^{k + 1}, u^{k + 1}\big) \leftarrow \big(P + \alpha F \big)^{-1} P(X^k, u^k) + \varepsilon^k,
	\end{aligned}\label{eq_approx_fixed_point}
\end{equation}
where $F$ and $P$ are the ones defined in (\ref{eq_F}) and (\ref{eq_P}), respectively, and $\varepsilon$ is an error component. In the literature, this methodology can be found under the name of \textit{inexact solves} or \textit{approximate proximal point} \cite{rockafellar1976monotone}. In the work of Eckstein and Bertsekas \cite{eckstein1992douglas}, they have shown that the approximate scheme (\ref{eq_approx_fixed_point}) converges as long as the error component is summable, i.e.
\begin{equation}
	\begin{aligned}
 	& \sum_{k=1}^{\infty} \left\Vert \varepsilon^k \right\Vert_2 < \infty.
	\end{aligned}\label{eq_summable}
\end{equation}

In the context of the approximate projection onto the p.s.d. cone, condition (\ref{eq_summable}) can be expressed in terms of the smallest eigenvalue of the truncated decomposition for each iterate. Let $\lambda_r^k$ denote the smallest eigenvalue computed at the $k^{\text{th}}$ iteration of the algorithm (\ref{eq_approx_fixed_point}), which corresponds to the $r^{\text{th}}$ largest eigenvalue at that iteration. Analogous to \cite{eckstein1992douglas}, given a target-rank $r$, the fixed point iteration (\ref{eq_approx_fixed_point}) will converge to a fixed point as long as
\begin{equation}
	\begin{aligned}
 	& (n - r) \sum_{k=1}^{\infty} \left\Vert \max\{\lambda_r^k, 0\} \right\Vert_2 < \infty
	\end{aligned}\label{eq_conv_condition}
\end{equation}
and a fixed point exists. It is easily verifiable that for an arbitrarily fixed target-rank $r$, the iteration (\ref{eq_approx_fixed_point}) will never converge. For instance, if one fixes the target-rank to a value smaller than the rank of the optimal solution, the error component will remain above a threshold and the sequence of errors will not be summable. We refer to the target-rank as \textit{sufficient} if it satisfies the condition (\ref{eq_conv_condition}).

Since the minimal sufficient target-rank is not known \textit{a priori}, it is necessary to use an update mechanism that can guarantee the convergence of (\ref{eq_approx_fixed_point}). The strategy adopted in this paper starts the algorithm with a small target-rank and increase its value whenever necessary. The combined residual (\ref{eq_comb_residual}) will be used to describe the \textit{state} of the algorithm and to trigger the update of the target-rank. Given an initial target-rank, the sum of the subsequent combined residuals can either converge or diverge. Even if the sequence converges, it will not necessarily be monotonic. In this regard, instead of checking the convergence of the sequential iterates we are going to evaluate the residuals (\ref{eq_comb_residual}) within a window of size $\ell$.


In case the combined residuals converge according to a given tolerance, we need to examine the approximation error (\ref{eq_error_approx}). It follows from (\ref{eq_error_approx}) that if the approximation error is zero, the smallest eigenvalue of $X^*$ is less than or equal to zero and the truncated projection is no longer an approximation. Therefore, the inexact iteration (\ref{eq_approx_fixed_point}) has also converged to a fixed point of (\ref{eq_fixed_point}) and consequently an optimal primal-dual solution for the SDP problem of interest has been found. 


If the combined residual has converged with a target rank $r$ but the approximation error is greater than the tolerance, the target-rank needs to be updated. In this case, it is interesting to notice that even though the current iterate pair is not an optimal point, it does give a feasible primal-dual solution under the assumption of strong duality. By characterizing a fixed point of (\ref{eq_error_approx}), the current iterate will satisfy the linear constraints of the original SDP problem. Additionally, as was previously pointed out, the truncated projection maps onto the positive semidefinite cone. Therefore, given a target rank $r$, if a fixed point of (\ref{eq_error_approx}) is found and strong duality holds, one has a feasible point designated by $(X_{[r]}, y_{[r]})$.

The last possible case occurs when the combined residuals either stay stationary or diverge within the last $\ell$ iterates. In this case, the target-rank also needs to be updated. After updating the target-rank,  the process is repeated. The combination of \texttt{PD-SDP} and the target rank updating scheme is described in Algorithm 2 and will be referred to as \texttt{LR-PD-SDP}. In the worst case scenario, the target-rank will be updated until $r$ equals $n$ and the subsequent iterations of the algorithm will be equivalent to the ones in \texttt{PD-SDP}. Consequently, in this setting, \texttt{LR-PD-SDP} will converge to a fixed point of (\ref{eq_F}) under the same conditions of \texttt{PD-SDP}.
\begin{algorithm}[H]
	\caption{\texttt{\texttt{LR-PD-SDP}}}
	\label{algorithm2} 
	\begin{algorithmic}
		\STATE{\textbf{Given:} $\mathcal{M}$, $b \in \mathbb{R}^p$, $h \in \mathbb{R}^q$, $C \in \mathbb{S}^n$ and $r=1$. \vspace{0.2cm}}
        \WHILE{$(n - r) \lambda_r > \varepsilon_{\lambda}$ \vspace{0.2cm}}
	    \WHILE{$\epsilon_{\text{comb}}^k > \epsilon_{\text{tol}}$ \textbf{and} $\epsilon_{\text{comb}}^{k} < \epsilon_{\text{comb}}^{k - \ell}$  \vspace{0.2cm}}
        \STATE{$X^{k + \smath{1}} \hspace{0.1cm} \leftarrow \textbf{aproj}_{\mathbb{S}_+^n} (X^k - \alpha (\mathcal{M}^*(y^k) + C), \hspace{0.05cm} r)$ \vspace{0.2cm}} \COMMENT{Approximate primal step}
        \STATE{$y^{k + \smath{1/2}} \hspace{-0.05cm} \leftarrow y^{k} + \alpha \mathcal{M}(2 X^{k + \smath{1}} - X^k)$ \vspace{0.1cm}}\COMMENT{Dual step part 1}
        \STATE{$y^{k + \smath{1}} \hspace{0.2cm} \leftarrow y^{k + \smath{1/2}} - \alpha \hspace{0.05cm} \textbf{proj}_{{=b \atop \leq h}}(y^{k + \smath{1/2}} / \alpha)$ }\COMMENT{Dual step part 2 \vspace{0.2cm}}
        \ENDWHILE{\vspace{0.1cm}}
        \IF{$\epsilon_{\text{comb}}^k < \epsilon_{\text{tol}}$}
        \STATE{$(X_{[r]}, y_{[r]}) \leftarrow (X^{k + 1}, y^{k + 1})$}\COMMENT{Save feasible solution}
        \ENDIF{}
        \STATE{$r \leftarrow 2 \hspace{0.05cm} r $ \vspace{0.1cm}}\COMMENT{\textit{Target-rank} update}
        \ENDWHILE{\vspace{0.2cm}}
		\RETURN{$(X^{k + 1}, y^{k + 1})$}
	\end{algorithmic}
\end{algorithm}

Each iteration of \texttt{LR-PD-SDP} has a computational complexity of $\mathcal{O}(n^2 r)$ as opposed to $\mathcal{O}(n^3)$ achieved by \texttt{PD-SDP}. Additionally, if one doubles the target-rank whenever necessary, the updating procedure can be carried out $\mathcal{O}(\text{log}(n))$ times. Usually, \texttt{LR-PD-SDP} will require more iterations to reach convergence than \texttt{PD-SDP}. On the other hand, \texttt{LR-PD-SDP} induces the rank of the iterates $X^k$ to remain small. As it is illustrated in Figure 2, \texttt{LR-PD-SDP} avoids the presence of high rank iterates, as happens with \texttt{PD-SDP}. Consequently, if the problem of interest has a low-rank solution, the \texttt{LR-PD-SDP} will terminate much faster than \texttt{PD-SDP}.
\begin{figure}[H]
\centering
\includegraphics[width=0.7\linewidth]{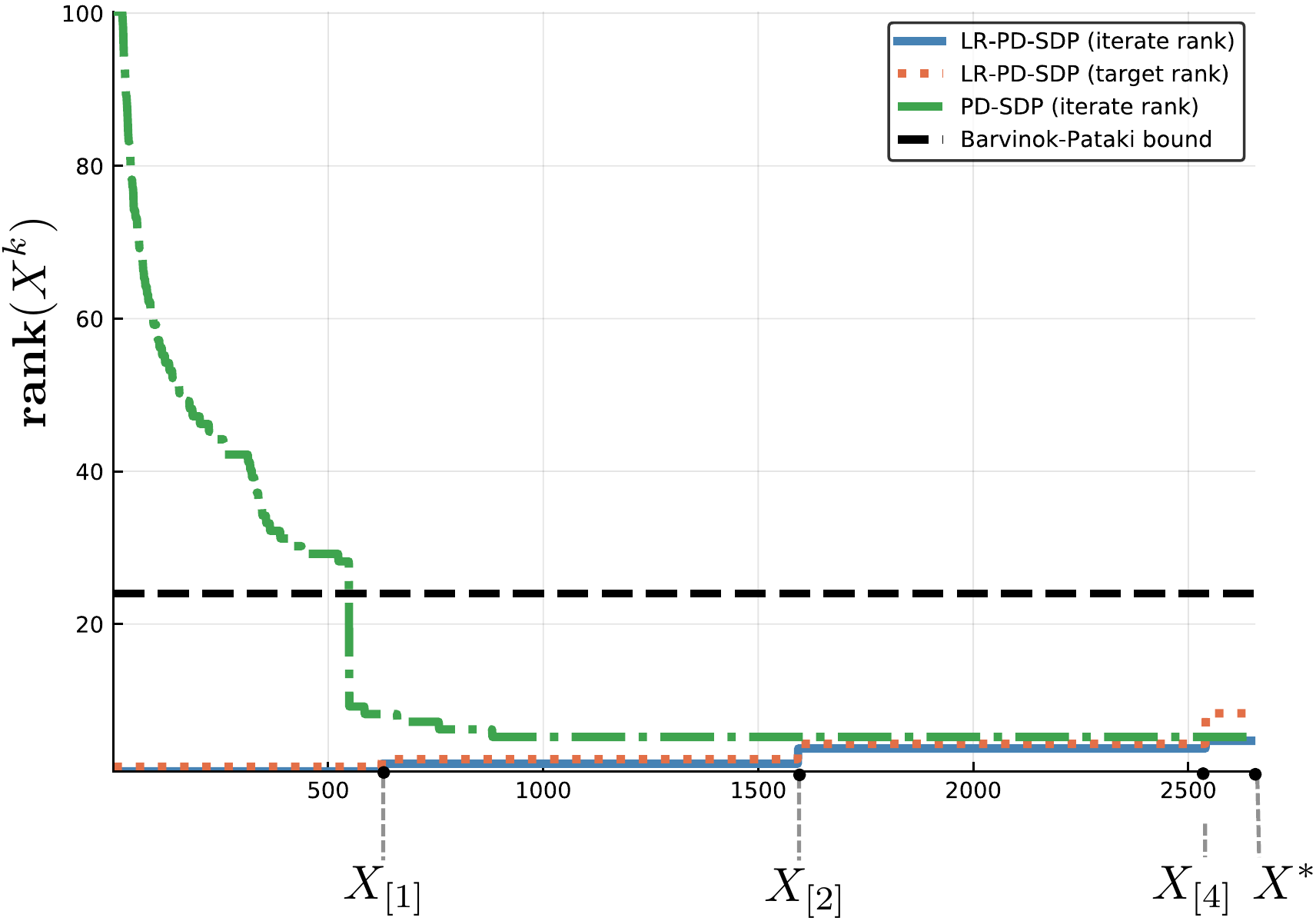}
\caption{Comparison of the rank path of the iterates $X^k$ for both \texttt{PD-SDP} and \texttt{LR-PD-SDP} methods. Additionally, the sequence of primal intermediate feasible solution found by \texttt{LR-PD-SDP} are represented as the points $X_{[1]}, X_{[2]}$ and $X_{[4]}$.}
\end{figure}




\section{\texttt{ProxSDP} solver}

The complete implementation of the \texttt{LR-PD-SDP} algorithm is available online at

\vspace{0.3cm}
\centerline{\href{https://github.com/mariohsouto/ProxSDP.jl}{\texttt{https://github.com/mariohsouto/ProxSDP.jl}}}
\vspace{0.3cm}
\noindent This project includes usage examples and all the data and scripts needed for next sections' benchmarks. The solver was completely written in the Julia language \cite{bezanson2017julia}, making extensive usage of its linear algebra capabilities. The use of sparse matrix operations were crucial to achieve good performance on manipulations involving the linear constraints. Additionally, dense linear algebra routines relying on BLAS \cite{lawson1979basic} were heavily used, just like multiple in-place operation to avoid unnecessary memory allocations. The built-in wrappers over LAPACK \cite{anderson1990lapack} and BLAS made it very easy to write high performance code. In particular, the ARPACK wrapper, used to efficiently compute the largest eigenvalues, was modified to maximize in-place operations and avoid unnecessary allocations.

Instead of writing a solver interface from scratch, we used the package MathOptInterface.jl (MOI) that abstracts solver interfaces. In doing that, we were able to write problems only once and test them in all available solvers. Moreover, having a MOI based interface means that the \texttt{ProxSDP} solver is available through the modeling language JuMP \cite{dunning2017jump}.

\section{Case studies}

In this section, we will present three SDP problems to serve as background for comparison between SDP solvers. The main goal of these experiments is to show how \texttt{LR-PD-SDP} outperforms the state-of-the-art solvers in the low-rank setting. In this sense, the numerical experiments are focused on semidefinite relaxation problems. On these SDP relaxations, the original problem one is interested in solving is nonconvex and it can be formulated as an SDP plus a rank constraint of the form $\textbf{rank}(X)=d$, where $d$ usually assumes a small value. Unfortunately, the rank constraint makes the problem extremely hard to solve and any exact algorithm has doubly exponential complexity \cite{chistov1984complexity}. The SDR avoids this problem by simply dropping the rank constraint and solving the remaining problem via semidefinite programming. Usually SDRs admit low-rank solutions, even without the presence of the rank constraint, making this the ideal case for testing the \texttt{LR-PD-SDP} algorithm.

In the presented experiments, we are going to consider a default numerical tolerance of $\epsilon_{\text{tol}}=10^{-3}$. As any first-order method, both \texttt{PD-SDP} and \texttt{LR-PD-SDP} may require a large number of iterations to converge to a higher accuracy \cite{eckstein1998operator, boyd2011distributed}. All the following tests were made using a Intel(R) Core(TM) i7-5820K CPU 3.30GHz (12 cores) Linux workstation with 62 Gb of RAM. In the following benchmarks for \texttt{PD-SDP} and \texttt{LR-PD-SDP} the Julia version used was compiled with Intel's MKL. The maximum running time for all experiments was set to 1200s.

\subsection{Graph equipartition}

Consider the undirected graph $G = (V, E)$ where $V$ is the set of vertices, $E$ is the set of edges, $n$ is the total number of edges and a cut $(S, S')$ is a disjoint partition of $V$. Let $x \in \{-1, +1\}^n$ such that 
\begin{equation*}
\begin{aligned}
& x_i = \left\{
\begin{array}{ll}
+1, \hspace{0.3cm} \text{if} \hspace{0.1cm} x_i \in S, \\
-1, \hspace{0.3cm} \text{if} \hspace{0.1cm} x_i \in S',
\end{array}
\right. \hspace{0.2cm} \forall \hspace{0.1cm} i=1, \cdots, n.
\end{aligned}
\end{equation*}
Given a set of weights $w$, the quantity $\tfrac{1}{4} \sum_{(i, j) \in E} w_{i j} (1 - x_i x_j)$ is called the weight of the cut $(S, S')$. The \textit{graph equipartition} problem aims to find the cut with maximum weight on a given graph such that both partitions of the graph have the same cardinality. This problem can be formulated as the following combinatorial optimization problem
\begin{equation*}
\begin{aligned}
& \underset{x}{\text{maximize}}
& & \tfrac{1}{4} \sum_{(i, j) \in E} w_{i j} (1 - x_i x_j)\\
& \text{subject to}
& & x_i \in \{-1, +1\}, \hspace{0.4cm} \forall \hspace{0.1cm} i = 1,\cdots, n,\\
&&& \sum_{i=1}^n x_i=0.
\end{aligned}\label{eq_max_cut_binary}
\end{equation*}

The binary constraints $x \in \{-1, +1\}^n$, can be expressed as a nonconvex equality constraints of the form $x_i^2=1 \hspace{0.1cm} \forall \hspace{0.1cm} i = 1, \cdots, n$. By lifting the decision variables to the space of the symmetric matrices $X \in \mathbb{S}_+^n$ and introducing a rank one constraint, the graph equipartition problem can be formulated as follows
\begin{equation*}
\begin{aligned}
& \underset{X \in \mathbb{S}_+^n}{\text{minimize}}
& & \textbf{tr}(W X)\\
& \text{subject to}
& & \textbf{tr}(\mathds{1}_{n \times n} X) = 0, \\
&&& \textbf{diag}(X) = 1,\\
&&& X \succeq 0,\\
&&& \textbf{rank}(X)=1,
\end{aligned}\label{eq_mimo_rank}
\end{equation*}
where the symmetric matrix $W$ is composed by the original weights $w$ and $\mathds{1}_{n \times n}$ denotes a $n \times n$ matrix filled with ones. By dropping the rank constraint, one obtains an SDP relaxation. For more details on graph partition problems, the reader should refer to \cite{karisch2000solving}.

\begin{table}[h!]
\centering
\begin{tabular}{|c|c|c|c|c|c|c|}
\hline
\textbf{n} & \textbf{sdplib} & \textbf{SCS} & \textbf{CSDP} & \textbf{MOSEK} & \textbf{PD-SDP} & \textbf{LR-PD-SDP} \\ \hline
124        & gpp124-1        & 1.6          & 0.4           & \textbf{0.2}   & 0.7             & 0.9                \\ \hline
124        & gpp124-2        & 1.5          & 0.4           & 0.3            & 0.5             & \textbf{0.2}       \\ \hline
124        & gpp124-3        & 1.6          & 0.3           & \textbf{0.2}   & 0.6             & \textbf{0.2}       \\ \hline
124        & gpp124-4        & 1.7          & 0.5           & 0.3            & 0.6             & \textbf{0.2}       \\ \hline
250        & gpp250-1        & 21.4         & 2.9           & \textbf{0.9}   & 3.7             & 1.4                \\ \hline
250        & gpp250-2        & 7.8          & 2.2           & \textbf{1.1}   & 4.1             & 1.2                \\ \hline
250        & gpp250-3        & 12.6         & 2.1           & \textbf{0.9}   & 3.4             & \textbf{0.9}       \\ \hline
250        & gpp250-4        & 16.4         & 2.2           & 0.9            & 3.8             & \textbf{0.6}       \\ \hline
500        & gpp500-1        & 134.2        & 59.1          & 8.2            & 22.7            & \textbf{5.6}       \\ \hline
500        & gpp500-2        & 97.4         & 12.2          & 8.6            & 21.5            & \textbf{6.1}       \\ \hline
500        & gpp500-3        & 64.4         & 12.1          & 8.9            & 15.5            & \textbf{4.4}       \\ \hline
500        & gpp500-4        & 71.4         & 13.4          & 8.7            & 15.4            & \textbf{6.5}       \\ \hline
801        & equalG11        & 324.2        & 47.3          & 32.4           & 84.3            & \textbf{11.3}      \\ \hline
1001       & equalG51        & 425.1        & 98.7          & 83.4           & 113.5           & \textbf{22.5}      \\ \hline
\end{tabular}
\caption{Comparison of running times (seconds) for the SDPLIB's graph equipartition problem instances.}
\end{table}

\textit{Problem instances:} \hspace{0.3cm} Graph equipartition instances from the SDPLIB \cite{borchers1999sdplib} problem set were used to evaluate the performance of the proposed methods. As Table 1 shows, for smaller instances Mosek solver is slightly faster. For larger instances such as equalG11 and equalG51, \texttt{LR-PD-SDP} outperforms all other considered methods with a considerable margin. Furthermore, without exploiting the low-rank structure of the problem, \texttt{PD-SDP} fails to scale as the number of edges increases. 

\newpage

\subsection{Sensor network localization}

Now consider the problem of estimating the position of a set of sensors on a $d$-dimensional plane \cite{alfakih1999solving}. Let $a_1, \dots, a_m \in \mathbb{R}^d$ be a set of anchor points in which the positions are known and let $x_1, \dots, x_n \in \mathbb{R}^d$ be a set of sensor points that have unknown positions. Given an incomplete set of Euclidean distances between sensors and between sensors and anchors, the goal is to find the true positions of each sensor. This problem, known as \textit{sensor network localization}, is originally formulated as the following quadratic constrained program:
\begin{equation}
\begin{aligned}
& \underset{x_1, \cdots, x_n \in \mathbb{R}^{d}}{\text{find}}
& & x_1, \cdots, x_n\\
& \text{subject to}
& & \left\Vert x_i - x_j \right\Vert_2^2 = w_{ij}^2, \hspace{0.3cm} \forall \hspace{0.1cm} (i, j) \in \Omega_s,\\
&&& \left\Vert a_k - x_j \right\Vert_2^2 = \tilde{w}_{kj}^2, \hspace{0.3cm} \forall \hspace{0.1cm} (k, j) \in \Omega_a,\\
\end{aligned}
\label{ed_quad_const}
\end{equation}
where the distances between sensor $i$ and sensor $j$ is denoted by $w_{ij}$ and the distance between anchor $k$ and sensor $j$ is denoted by $\tilde{w}_{kj}$. The indexes of the distances that are known are either in the set $\Omega_s$ or in the set $\Omega_a$. Unfortunately, solving (\ref{ed_quad_const}) is NP-hard \cite{boyd2004convex}.

We can formulate the nonconvex problem (\ref{ed_quad_const}) as a rank constrained semidefinite problem \cite{so2007theory}. In order to start building this alternative formulation, consider the matrices $X \in \mathbb{R}^{d \times n}$ and $Y \in \mathbb{S}^n$ as
\begin{equation*}
\begin{aligned}
& X = 
\left[
  \begin{array}{ccc}
    \vertbar &         & \vertbar \\
    x_{1}    & \cdots & x_{n}    \\
    \vertbar &       & \vertbar 
  \end{array}
\right]  \hspace{0.1cm} \text{and} \hspace{0.1cm}
Y = X^T X = \begin{bmatrix}
    x_{1}^T x_{1} & x_{1}^T x_2  & \dots  & x_{1}^T x_n \\
    x_{2}^T x_1 & x_{2}^T x_{2}  & \dots  & x_{2}^T x_n \\
    \vdots & \vdots  & \ddots & \vdots \\
    x_{n}^T x_1 & x_{n}^T x_2  & \dots  & x_{n}^T x_{n}
\end{bmatrix} .
\end{aligned}
\end{equation*}
Now let $E^{(i, j)} \in \mathbb{S}^n$ be filled with zeros except for the following entries: $E^{(i, j)}_{i,i}=1$, $E^{(i, j)}_{j,j}=1$, $E^{(i, j)}_{i,j}=-1$ and $E^{(i, j)}_{j,i}=-1$. With this setting, the constraints that represent the distance between sensors can be formulated as
\begin{equation}
\begin{aligned}
& \textbf{tr}(E^{(i, j)} Y) = \omega_{ij}^2, \hspace{0.4cm} \forall \hspace{0.1cm} (i, j) \in \Omega_s.
\end{aligned}\label{eq_dist_sensors}
\end{equation}
Similarly, let $Z \in \mathbb{S}^{d + n}$ be the matrix
\begin{equation*}
\begin{aligned}
& Z = \begin{bmatrix}
    I_{d \times d} & X \\
    X^T & Y
  \end{bmatrix}.
\end{aligned}
\end{equation*}
Additionally, let $U^{(k, j)} \in \mathbb{S}^{d + n}$ be filled with zeros except for the entries: $U^{(k, j)}_{1,1}=a_k^T a_k$, $U^{(k, j)}_{d+j,d +j}=1$, $U^{(k, j)}_{1:d,d + j}=- a_k$ and $U^{(k, j)}_{d + j, 1:d}=- a_k^T$. The constraints regarding the distances between sensors and anchors can be formulated as
\begin{equation}
\begin{aligned}
& \textbf{tr}(U^{(k, j)} Z) = \tilde{\omega}_{kj}^2, \hspace{0.4cm} \forall \hspace{0.1cm} (k, j) \in \Omega_a.
\end{aligned}\label{eq_dist_anchors}
\end{equation}


Using (\ref{eq_dist_sensors}), (\ref{eq_dist_anchors}) and the Schur complement of $Y \succeq X^T X$ \cite{so2007theory}, the network localization problem can be formulated as
\begin{equation}
\begin{aligned}
& \underset{Z \in \mathbb{S}_+^{d + n}}{\text{find}}
& & Z\\
& \text{subject to}
& & Z = \begin{bmatrix}
    I_{d \times d} & X \\
    X^T & Y
  \end{bmatrix} \succeq 0,\\
&&& \textbf{tr}(E^{(i, j)} Y) = \omega_{ij}^2, \hspace{0.4cm} \forall \hspace{0.1cm} (i, j) \in \Omega_s,\\
&&& \textbf{tr}(U^{(k, j)} Z) = \tilde{\omega}_{kj}^2, \hspace{0.4cm} \forall \hspace{0.1cm} (k, j) \in \Omega_a,\\
&&& \textbf{rank}(Y) = d.
\end{aligned}
\label{eq_network_loc}
\end{equation}
If a unique solution for the given set of distances exists, the SDR obtained by dropping the rank constraint in (\ref{eq_network_loc}) will be exact \cite{so2007theory}.

\textit{Problem instances:} \hspace{0.3cm} In a set of numerical simulation, we randomly generate anchor points and distances measurements. Each anchor and sensor has its position in the two-dimensional Euclidean plane, i.e. $d=2$. In this sense, if the relaxation is exact the optimal solution $Y^*$ must have a rank of two. This property justifies the \texttt{LR-PD-SDP} outperforming other solvers when number of sensors grow, as it can be seen in Table 2. The importance of exploiting the low rank structure can be verified by observing that \texttt{PD-SDP} does not efficiently scale as $n$ increases.

\begin{table}[h!]
\centering
\begin{tabular}{|c|c|c|c|c|c|}
\hline
\textbf{n} & \textbf{SCS} & \textbf{CSDP} & \textbf{MOSEK} & \textbf{PD-SDP} & \textbf{LR-PD-SDP} \\ \hline
50         & 0.2          & 0.2           & \textbf{0.1}   & 0.5             & 0.6                \\ \hline
100        & \textbf{0.8} & 4.5           & 0.9            & 6.1             & 1.6                \\ \hline
150        & \textbf{2.6} & 28.1          & 3.2            & 14.4            & 3.6                \\ \hline
200        & 6.4 & 89.8          & 11.2           & 32.3            & \textbf{6.1}                \\ \hline
250        & 12.1         & 239.2         & 36.4           & 52.9            & \textbf{7.9}       \\ \hline
300        & 28.7         & timeout       & 85.2           & 96.6            & \textbf{13.5}      \\ \hline
\end{tabular}
\caption{Comparison of running times (seconds) for randomized network localization problem instances.}
\end{table}

\subsection{MIMO detection}

Consider an application in the field of wireless communication known in the literature as binary multiple-input multiple-output (MIMO) \cite{jalden2003semidefinite, jalden2008diversity}. As in several MIMO applications, one needs to send and receive multiple data signals over the same channel with the presence of additive noise. The binary MIMO can be modeled as:
\begin{equation*}
	\begin{aligned}
	& y = H x + \varepsilon,
	\end{aligned}
\end{equation*}
where $y \in \mathbb{R}^m$ is the received signal, $H \in \mathbb{R}^{m \times n}$ is the channel and $\varepsilon \in \mathbb{R}^m$ is an i.i.d. Gaussian noise with variance $\sigma^2$. The signal, which is unknown to the receiver, is represented by $x \in \{-1, +1\}^n$.

Assuming the noise distribution is known to be Gaussian, a natural approach is to compute the maximum likelihood estimate of the signal by solving the optimization problem:
\begin{equation}
\begin{aligned}
& \underset{x}{\text{minimize}}
& & \left\Vert H x - y \right\Vert_2^2\\
& \text{subject to}
& & x \in \{-1, +1\}^n.
\end{aligned}\label{eq_binary_least_squares}
\end{equation}
At first sight, the structure of the problem is very similar to a standard least squares problem. However, the unknown signal is constrained to be binary, which makes the problem nonconvex and dramatically changes the problem's complexity. More precisely, solving \eqref{eq_binary_least_squares} is known to be NP-hard \cite{verdu1989computational}.

By using a similar technique as in the graph equipartition problem, one can reformulate \eqref{eq_binary_least_squares} as the following rank constrained semidefinite problem:
\begin{equation*}
\begin{aligned}
& \underset{X \in \mathbb{S}_+^{n+1}}{\text{minimize}}
& & \textbf{tr}(L X)\\
& \text{subject to}
& & \textbf{diag}(X) = 1,\\
&&& X_{n+1, n+1} = 1,\\
&&& -1 \leq X \leq 1,\\
&&& X \succeq 0,\\
&&& \textbf{rank}(X) = 1.
\end{aligned}
\end{equation*}
where the decision variable $X$ is a $n + 1 \times n + 1$ symmetric matrix and $L$ is given by:
\begin{equation*}
\begin{aligned}
L=
  \begin{bmatrix}
    H^{*}H & -H^{*}y \\
    -y^{*} H & y^{*} y
  \end{bmatrix}
\end{aligned}\label{eq_matrix_completion}
\end{equation*}
Given the optimal solution $X^*$ for the relaxation, the solution for the original binary MIMO is obtained by slicing the last column as $x^* = X^*_{1:n, n + 1} $. For this particular problem, the SDR is known to be exact if the signal to noise ratio, $\sigma^{-1}$, is sufficiently large \cite{jalden2003semidefinite}. This implies that the rank of the optimal solution $X^*$ is guaranteed to be equal to one even without the rank constraint. This low-rank structure makes the ideal case study for the techniques proposed in this paper.

\textit{Problem instances:} \hspace{0.3cm} In order to measure the performance of the different methods, problem instances with large signal to noise ratio were randomly generated. For each instance, the channel matrix $H$ is designed as a $n \times n$ matrix with i.i.d. standardized Gaussian entries. The true signal $x^*$ was drawn from a discrete uniform distribution. Since a high signal to noise ratio was used to build the instances, all recovered optimal solutions are rank one solutions. In this setting, as it is illustrated in Table 3, \texttt{LR-PD-SDP} outperforms all other methods as the signal length increases. More surprisingly, \texttt{LR-PD-SDP} was able to solve large scale instances with $5000 \times 5000$ p.s.d. matrices. The bottleneck found while trying to optimize even larger instances was the amount of memory required by the \texttt{ProxSDP} solver. 

\begin{table}[h!]
\centering
\begin{tabular}{|c|c|c|c|c|c|}
\hline
\textbf{n} & \textbf{SCS} & \textbf{CSDP} & \textbf{MOSEK} & \textbf{PD-SDP} & \textbf{LR-PD-SDP} \\ \hline
100        & 1.5          & 1.2            & \textbf{0.1}   & \textbf{0.1}    & \textbf{0.1}       \\ \hline
500        & 277.8        & 27.4           & 2.3            & 3.1             & \textbf{1.1}       \\ \hline
1000       & timeout      & 97.2           & 15.6           & 16.5            & \textbf{4.7}       \\ \hline
2000       & timeout      & 473.6          & 117.5          & 115.9           & \textbf{38.9}      \\ \hline
3000       & timeout      & timeout        & 418.2          & 350.6           & \textbf{122.1}     \\ \hline
4000       & timeout      & timeout        & 976.8          & 906.5           & \textbf{258.3}     \\ \hline
5000       & timeout      & timeout        & timeout        & timeout         & \textbf{472.4}     \\ \hline
\end{tabular}
\caption{Running times (seconds) for MIMO detection with high SNR.}
\end{table}
 
\section{Conclusions and future work}

As a concluding remark, this work has proposed a novel primal-dual method that can efficiently exploit the low-rank structure of semidefinite programming problems. As it was illustrated by the case studies, the proposed technique can achieve up to one order of magnitude faster solving times in comparison to existing algorithms. Additionally, an open source solver, \texttt{ProxSDP}, for general SDP problems was made available. We hope that the results and tools contemplated in this work foster the use of semidefinite programming on new applications and fields of study.

One aspect of the proposed methodology not fully explored in this paper, is the value of the intermediate solutions found by \texttt{LR-PD-SDP}. For several applications, a suboptimal feasible solution may be useful. Particularly if one is interested in solving a semidefinite relaxation, a suboptimal solution can be almost as useful as the optimal solution, with the advantage of requiring less computing time to be discovered. For instance, a branch-and-bound search method can benefit from lower bounds that a feasible semidefinite relaxation provides \cite{dong2016relaxing}. This ability of quickly generating high quality lower bounds via intermediate feasible solutions can enhance the already well known SDP property of approximate hard problems.

Another promising future line of work is the combination of chordal decomposition methods with the low-rank approximation presented in this work. If successful, this match would allow the exploitation of both sparsity and low-rank structure simultaneously.



\section*{Acknowledgments}

Firstly, we would like to thank the Brazilian agencies CNPq and CAPES for financial
support. We extend many thanks to all members of LAMPS (Laboratory of Applied Mathematical Programming and Statistics), in special Thuener Silva and Raphael Saavedra, for the daily support and fruitful discussions. We would also like to thank the developers of MOI and JuMP for making comparisons between solvers much easier and specially Benoît Legat for helping with \texttt{ProxSDP}'s MOI interface.

\bibliographystyle{IEEEtran}
\bibliography{Ref}
\vfill

\end{document}